\documentclass[final]{amsart}
\usepackage{url}

\newcommand{\eqinlaw}{\stackrel{\mathcal{L}}{=}}
\newcommand{\toinlaw}{\xrightarrow{\mathcal{L}}}
\newcommand{\E}[1]{\mathbf{E}\,{#1}}
\newcommand{\mun}[1]{\mu_n^{[{#1}]}}
\newcommand{\muz}[1]{\mu^{[{#1}]}(z)}
\newcommand{\sums}{\sum_{\substack{s_1+s_2+s_3=s\\ s_2, s_3 < s}}}

\newcommand{\A}{\mathcal{A}}
\newcommand{\N}{\mathcal{N}}
\newcommand{\wX}{\widetilde{X}}

\DeclareMathOperator{\Li}{Li}

\newtheorem{theorem}{Theorem}[section]
\newtheorem{proposition}[theorem]{Proposition}
\theoremstyle{remark}
\newtheorem{remark}[theorem]{Remark}

\numberwithin{equation}{section}

\begin{document}

\title{Destruction of very simple trees}
\author{James Allen Fill}
\author{Nevin Kapur}
\author{Alois Panholzer}
\date{August~4, 2005}

\email[James Allen Fill]{jimfill@jhu.edu}
\email[Nevin Kapur]{nkapur@cs.caltech.edu}
\email[Alois Panholzer]{Alois.Panholzer@tuwien.ac.at}

\urladdr[James Allen Fill]{http://www.ams.jhu.edu/\~{}fill/}
\urladdr[Nevin Kapur]{http://www.cs.caltech.edu/\~{}nkapur/}
\urladdr[Alois Panholzer]{http://info.tuwien.ac.at/panholzer/}

\thanks{The research of James Allen Fill was supported by \textsc{NSF}
  Grants \textsc{DMS}--0104167 and \textsc{DMS}--0406104, and by The
  Johns Hopkins University's Acheson J.~Duncan Fund for the
  Advancement of Research in Statistics.  Nevin Kapur's research was
  supported by \textsc{NSF} grant~0049092 and the Center for
  Mathematics of Information at the California Institute of Technology.}

\address[James Allen Fill]{Applied Mathematics and Statistics, The
  Johns Hopkins University, 3400 N.~Charles St., Baltimore MD
  21218-2682}
\address[Nevin Kapur]{Computer Science, California Institute of
  Technology, MC 256-80, 1200 E.~California Blvd., Pasadena CA 91125}
\address[Alois Panholzer]{Institut f{\"u}r Diskrete Mathematik und
  Geometrie, Technische Universit{\"a}t Wien, Wiedner Hauptstra{\ss}e
  8--10/104, A-1040 Wien, Austria}



\begin{abstract}
  We consider the total cost of cutting down a random rooted tree
  chosen from a family of so-called \emph{very simple trees} (which include
  ordered trees, \mbox{$d$-ary} trees, and Cayley trees); these form a
  subfamily of simply generated trees.  At each stage of the process an
  edge is chose at random from the tree and cut, separating the tree
  into two components.  In the one-sided variant of the process the
  component not containing the root is discarded, whereas in
  the two-sided variant both components are kept.  The process ends
  when no edges remain for cutting.  The cost of cutting an
  edge from a tree of size~$n$ is assumed to be~$n^\alpha$.  Using
  singularity analysis and the method of moments, we derive the
  limiting distribution of the total cost accrued in both variants of
  this process.  A salient feature of the limiting distributions
  obtained (after normalizing in a family-specific manner) is that
  they only depend  on~$\alpha$.
\end{abstract}

\maketitle


\section{Introduction}
\label{sec:introduction}

Consider the following process on a rooted tree with $n$~vertices.
Pick an edge uniformly at random and ``cut'' it, separating the tree
into a pair of rooted trees; the tree containing the root of the
original tree retains its root while the tree not containing the root
of the original tree is rooted at the vertex adjacent to the edge that
was cut.  In the \emph{one-sided} variant of the problem the tree not
containing the original root is discarded and the process is continued
recursively until the original root is isolated.  In the
\emph{two-sided} variant the process is continued recursively on each
of the rooted trees.  Assume that the cost incurred for selecting an
edge and splitting the tree is $t_n$.  In this paper we derive the
limiting distribution of the total cost accrued when the tree is a
random \emph{very simple tree} (defined below) and $t_n = n^\alpha$
for fixed $\alpha \geq 0$, for both the two-sided variant
(Theorems~\ref{thm:dist12},~\ref{thm:lt12},and~\ref{thm:moments-half})
  and the one-sided variant (Theorem~\ref{theo1}).  A salient feature
  of the limiting distributions obtained (after normalizing in a
  family-specific manner) is that they only depend on~$\alpha$.

In the one-sided variant, the case $t_n \equiv 1$ (i.e., $\alpha = 0$)
corresponds to the number of cuts required to disconnect the tree.
For this random variable, Meir and Moon~\cite{MR44:1598} derived the
mean and variance for Cayley trees; Chassaing and
Marchand~\cite{chassaing_marchand_cuts} derived the limiting
distribution for Cayley trees.  Panholzer obtained limiting
distributions for non-crossing trees~\cite{MR2042393} and very simple
families of trees~\cite{panholzer:2003}.  Recently Janson extended
these results to all simply generated families~\cite{janson-167}.

The interest in the two-sided variant stems from the fact that when
the very simple family is Cayley trees, the process is equivalent to
a probabilistic model (the ``random spanning tree model'')
involved in the Union--Find (or equivalence-finding) algorithm.  Knuth
and Sch{\"o}nhage~\cite{MR81a:68049} derived the expected value of the
cost in the cases (among others) $t_n \sim a \sqrt{n}$ and $t_n =
n/2$.  These results were later extended~\cite{MR89i:05024} to the
cases $t_n = n^\alpha$ when $\alpha > 1/2$ and $t_n = O(n^{\alpha})$
when $\alpha < 1/2$.  (Some of these expected values were rederived
using singularity analysis in~\cite{FFK}.)  In~\cite{math.PR/0406094},
Chassaing and Marchand derive limit laws for the costs considered by
Knuth and Sch{\"o}nhage.

We treat both variants of the destruction process using singularity
analysis~\cite{MR90m:05012}, a complex-analytic technique that relates
asymptotics of sequences to singularities of their generating
functions.  We rely on applicability of singularity analysis to the
Hadamard product (the term-by-term product) of sequences~\cite{FFK}
and the amenability of the generalized polylogarithm to singularity
analysis~\cite{MR2000a:05015}.

The organization of the paper is as follows.  In
Section~\ref{sec:very-simple-trees} we define families of very simple
trees, noting the key ``randomness-preservation'' property that is
crucial for the application of our methods.
Section~\ref{sec:preliminaries} establishes notation and other
preliminaries that will be used in the subsequent proofs.  In
Section~\ref{sec:toll-nalpha} the two-sided variant is considered, and
Section~\ref{sec:one-sided-destr} deals with the one-sided variant.

\medskip
\noindent
\emph{Notation.}  In the sequel we will use~$\ln$ to denote natural
logarithms and~$\log$ when the base of the logarithm does not matter.

\section{Very simple trees}
\label{sec:very-simple-trees}
An \emph{ordered tree} is a rooted tree in which the order of the
subtrees of each given node is relevant.  Given a sequence
$(\phi_i)_{i \geq 0}$ of nonnegative numbers (called a degree
generating sequence) with $\phi_0 = 1$, a \emph{simply generated
  family}~$\mathcal{F}$ of trees is obtained by assigning each ordered
tree $T$ the weight
\begin{equation*}
  w(T) := \prod_{v \in T} \phi_{d(v)},
\end{equation*}
where $d(v)$ is the outdegree of the node $v$.  Let $\mathcal{F}_n$
denote the set of trees in $\mathcal{F}$ with~$n$ nodes, and let $T_n$
denote the weighted number of trees in $\mathcal{F}_n$, i.e.,
\begin{equation*}
  T_n := \sum_{T \in \mathcal{F}_n} w(T).
\end{equation*}

A \emph{random simply generated tree} of size~$n$ is obtained by
assigning probability~$w(T)/T_n$ to the tree $T \in \mathcal{F}_n$.
Many combinatorially interesting families such as (unweighted) ordered
trees, Cayley trees, Motzkin trees, and $d$-ary trees are simply
generated.  It is also well known that simply generated trees correspond to
certain conditioned Galton--Watson trees; see the introductory section
of~\cite{janson-167} for the precise connection.  It is well-known
that the
generating function $T(z) := \sum_{n \geq 1} T_n z^n$
satisfies the functional equation
\begin{equation*}
   T(z) = z \Phi(T(z)),
\end{equation*}
where $\Phi(t) := \sum_{k \geq 0} \phi_k t^k$ is the degree generating
function of the family.  For further background on simply generated
trees we refer the reader to~\cite{MR80k:05043}.

In this paper we consider the subclass of simply generated families,
called \emph{very simple families}, that, among simply generated
families, are characterized by the following property.

\medskip
\begin{center}
\parbox[c]{.9\linewidth}{%
Choose a random simply generated tree from the family~$\mathcal{F}_n$
and then one of its $n-1$ edges uniformly at random.  Cutting this
edge produces a pair of trees of size~$k$ (the one that contains the
root) and $n-k$, as described in Section~\ref{sec:introduction}.  Then
the subtrees themselves are random simply generated trees from the
family~$\mathcal{F}_k$ and $\mathcal{F}_{n-k}$.
}
\end{center}

\medskip

It is clear that the ``randomness-preservation'' property of very
simple trees allows a simple recursive formulation [see~\eqref{eq:1}
and~\eqref{eqno1}] of the total cost of destroying such a tree.

Panholzer~\cite[Lemma~1]{panholzer:2003} characterized the degree
generating functions of very simple trees; the relevant constraints
are summarized in Table~\ref{tab:sp}.  

\subsection{Singular expansions}
\label{sec:singular-expansions}

As is usual for treatment of
simply generated families, let $\tau$ denote the unique root of $t
\Phi'(t) = \Phi(t)$ with $0 < t < R$, where $R$ is the radius of
convergence of the series~$\Phi$.  Let $\rho := \tau/\Phi(\tau)$.  Let $Z :=
1-\rho^{-1}z$, and let $\A$ denote a generic power series in $Z$,
possibly different at each occurrence.  Then as $z \to \rho$, the
dominant singularity for~$T(z)$, a singular expansion for $T$
is~\cite[Theorem~VII.2]{flajolet:_analy_combin}
\begin{equation}
  \label{eq:7}
  T(z) \sim \tau - b\rho^{1/2}Z^{1/2} + Z\A + Z^{3/2}\A,
\end{equation}
where $b := \Phi(\tau)\sqrt{\frac{2}{\tau \Phi''(\tau)}}$.  immediately from singularity analysis that
\begin{equation}
  \label{eq:15.1}
  T_n \sim c \rho^{-n} n^{-3/2}(1 +  n^{-1}\N),
\end{equation}
where $c := b\rho^{1/2}/(2\sqrt\pi) =
[\Phi(\tau)/2\pi\Phi''(\tau)]^{1/2}$.  In the sequel we will also use
\begin{equation}
  \label{eq:15}
  \sigma^2 := \tau^2 \frac{\Phi''(\tau)}{\Phi(\tau)}.
\end{equation}

Differentiating the expansion~\eqref{eq:7}
term-by-term~\cite[Theorem~6]{FFK} we get
\begin{equation*}
  T'(z) \sim \frac{b}2 \rho^{-1/2}Z^{-1/2} + \A + Z^{1/2}\A.
\end{equation*}
Since $z = \rho - \rho Z$,
\begin{equation}
  \label{eq:8}
  z T'(z) \sim \frac{b}{2} \rho^{1/2} Z^{-1/2} + \A + Z^{1/2}\A.
\end{equation}
The constants~$a_0$ and~$a_1$ described by Table~\ref{tab:sp} are
fundamental constants for our analysis; see
especially~\eqref{eq:p_nk}.  Using~\eqref{eq:7} and~\eqref{eq:8} we
get
\begin{equation}
  \label{eq:9}
  1 + 2 a_0 T(z) + a_1 z T'(z) \sim a_1 \rho^{1/2} \frac{b}{2} Z^{-1/2} + \A
  + Z^{1/2}\A
\end{equation}
and
\begin{equation*}
  1 - a_1 T(z) \sim (1- a_1\tau) + a_1 b \rho^{1/2} Z^{1/2} + Z\A +
  Z^{3/2}\A.
\end{equation*}
It is easily verified that for each very simple family $1 - a_1\tau =
0$ (this fact will be used numerous times in  subsequent
calculations), so that the constant term vanishes in the singular
expansion of $1 - a_1 T(z)$.  This leads to
\begin{equation*}
  z[1 - a_1 T(z)] \sim \rho^{3/2} a_1 b Z^{1/2} + Z\A + Z^{3/2}\A
\end{equation*}
and consequently
\begin{equation}
  \label{eq:10}
  z^{-1}[1 - a_1 T(z)]^{-1} \sim \rho^{-3/2} a_1^{-1} b^{-1} Z^{-1/2}
  + \A + Z^{1/2}\A.
\end{equation}

We will also need
\begin{equation}
  \label{eq:27}
  \frac1{T(z)} \sim {\tau}^{-1} + \frac{b\rho^{1/2}}{\tau^2}
  Z^{{1}/{2}} + Z \A + Z^{{3}/{2}} \A,
\end{equation}
which follows from~\eqref{eq:7}.

\section{Preliminaries}
\label{sec:preliminaries}

Throughout, $\eqinlaw$ denotes equality in law (or distribution) and
$\toinlaw$ denotes convergence in law.  Recall that
the \emph{Hadamard product} of two power series $f$ and $g$, denoted
by $f \odot g$, is the power series defined by
\begin{equation*}
  ( f \odot g)(z) \equiv f(z) \odot g(z) := \sum_{n} f_n g_n z^n,
\end{equation*}
where
\begin{equation*}
 f(z) = \sum_{n} f_n z^n \quad \text{ and } \qquad g(z) = \sum_{n}
 g_n z^n.
\end{equation*}

\subsection{Two-sided destruction}
\label{sec:two-sided-destr}

The cost of cutting down a very simple tree of size~$n$, call it
$X_n$, satisfies the distributional recurrence
\begin{equation}
  \label{eq:1}
  X_n \eqinlaw X_{K_n} + X_{n-K_n}^* + t_n, \quad n \geq 2; \qquad X_1
  = t_1,
\end{equation}
where $t_n$, for $n \geq 2$, is the toll for cutting an edge from a
tree of size~$n$.  Here $K_n$, the (random) size of the tree
containing the root, is independent of $(X_j)_{j \geq 1}$ and
$(X_j^{*})_{j \geq 1}$, which are independent copies of each other.
The \emph{splitting probabilities} are given by
\begin{equation}
  \label{eq:p_nk}
  \Pr[ K_n = k ] =: p_{n,k} = (a_1 k + a_0) \frac{T_k T_{n-k}}{(n-1)
  T_n}, \quad k=1,\ldots,n-1.
\end{equation}
Table~\ref{tab:sp} gives the constants~$a_1$ and~$a_0$ for each type
of very simple family; see~(14)--(16) in~\cite{panholzer:2003}. Here
$\alpha_i := \phi_{i+1}/\phi_i$, $i=0,1$, where $(\phi_i)_{i
  \geq 0}$ is the degree generating sequence of the simply generated
tree.
\begin{table}[htbp]
  \centering
  \begin{tabular}[p]{ccccc}
    Family  & Generating function  & Constraints & $a_1$   & $a_0$ \\
    \hline
    A       & $e^{\alpha_0t}$ &  & $\alpha_0$    & 0     \\
    B       & $(1 + \frac{\alpha_0t}d)^d$ & $d \geq 2$ & $\alpha_0
    \frac{d-1}{d}$ &    $\frac{\alpha_0}d$       \\ 
    C       & $[1 - (2\alpha_1 -
    \alpha_0)t]^{-\frac{\alpha_0}{2\alpha_1-\alpha_0}}$ & $2\alpha_1 -
    \alpha_0 > 0$ & $2\alpha_1$  &  $-(2\alpha_1 - \alpha_0)$    \\
            &       &     &            &
  \end{tabular}
  \caption{Generating functions for very simple families.  For each
    family, $\alpha_0 > 0$ is also a constraint.}
  \label{tab:sp}
\end{table}
It is easy to check that family~A is Cayley trees, family~B is $d$-ary
trees, and family~C contains unweighted ordered trees.  (As it turns out, the
distributional recurrence for Cayley trees is identical to the one
obtained for the Union--Find process studied
in~\cite{MR81a:68049,MR89i:05024,FFK}---see Remark~\ref{rem:symmetry}
below.)

Define $\mun{s} := \E{X_n^s}$.  Taking $s$th powers of both sides
of~\eqref{eq:1} and taking expectations by conditioning on~$K_n$, we get
\begin{equation}
  \label{eq:2}
  \mun{s} = \sum_{k=1}^{n-1} p_{n,k} ( \mu_k^{[s]} + \mu_{n-k}^{[s]} ) +
  r_n^{[s]}, \quad n \geq 2,
\end{equation}
where
\begin{equation}
  \label{eq:29}
  r_n^{[s]} := \sums \binom{s}{s_1,s_2,s_3} t_n^{s_1} \sum_{k=1}^{n-1}
  p_{n,k} \mu_{k}^{[s_2]} \mu_{n-k}^{[s_3]},
\end{equation}
and $\mu_1^{[s]} = t_1^s$.  Define generating functions
\begin{equation*}
  \muz{s} := \sum_{n \geq 1} \mun{s} T_n z^n, \qquad t(z) := \sum_{n \geq 1}
  t_n z^n, \qquad T(z) := \sum_{n \geq 1} T_n z^n.
\end{equation*}
[Observe that $\mu^{[0]}(z) = T(z)$.]
Multiply~\eqref{eq:2} by $(n-1)T_n z^n$ and sum over $n \geq 2$.
The resulting left side is
\begin{equation*}
  \sum_{n \geq 2} (n-1) T_n \mun{s} z^n = \sum_{n \geq 1} (n-1) T_n \mun{s}
  z^n = z \partial_z \muz{s} - \muz{s},
\end{equation*}
where $\partial_z$ denotes derivative with respect to~$z$.
Similarly, the resulting first term on the right side is
\begin{equation*}
  a_1\left[ z \left( \partial_z \muz{s} \right) T(z) + z T'(z)\muz{s}
  \right]  + 2 a_0 \muz{s} T(z).
\end{equation*}
The resulting second term on the right side is
\begin{align}
 & r^{[s]}(z) := \sum_{n \geq 2} (n-1) T_n r_n^{[s]} z^n = \sum_{n \geq 1}
  (n-1) T_n r_n^{[s]} z^n \notag \\
  \label{eq:3}
 & = \sums \binom{s}{s_1,s_2,s_3} t^{\odot s_1}(z) \odot \left[ a_1 z
  \left( \partial_z \muz{s_2} \right) \muz{s_3} + a_0 \muz{s_2}
  \muz{s_3} \right].
\end{align}
Thus~\eqref{eq:2} translates to
\begin{multline*}
  z \partial_z \muz{s} - \muz{s} \\ = a_1\left[ z \left( \partial_z
  \muz{s} \right) T(z) + z T'(z)\muz{s}  \right]  + 2 a_0 \muz{s} T(z)
  + r^{[s]}(z),
\end{multline*}
i.e.,
\begin{equation}
  \label{eq:3.1}
  \partial_z \muz{s} + p(z) \muz{s} = g^{[s]}(z),
\end{equation}
where
\begin{equation}
  \label{eq:4}
  p(z) := - \frac{1 + 2a_0 T(z) + a_1 z T'(z)}{z[1 - a_1T(z)]}
\end{equation}
and
\begin{equation}
  \label{eq:5}
  g^{[s]}(z) := \frac{r^{[s]}(z)}{z[1-a_1T(z)]},
\end{equation}
with $\mu^{[s]}(0) = 0$.  By variation of parameters (see, for
example,~\cite[2.1-(22) and
Problem~2.1.21]{boyce86:_elemen_differ_equat}, the general solution to
the first-order linear differential equation~\eqref{eq:3.1} is given by
\begin{equation}
  \label{eq:6.1}
  \muz{s} = A^{[s]}(z) \exp\left[ - \int_{z_0}^z p(t)\,dt \right],
\end{equation}
where
\begin{equation}
  \label{eq:6}
  A^{[s]}(z) := \int_0^z g^{[s]}(t) \exp\left[ \int_{z_0}^t p(u)\,du
  \right]\,dt + \beta_s,
\end{equation}
with $z_0$ chosen as follows and $\beta_s$ an arbitrary constant.

The integrand~$p(z)$ defined at~\eqref{eq:4} and appearing
in~\eqref{eq:6.1}--\eqref{eq:6} is asymptotic to~$-1/z$ as $z \to 0$
and has~[see~\eqref{eq:38} below] another singularity at~$z=\rho$.
In~\mbox{\eqref{eq:6.1}--\eqref{eq:6}} we may choose (and fix) $z_0$
arbitrarily from the punctured disc of radius~$\rho$ centered at the
origin.  Then, in~\eqref{eq:6}, as $t \to 0$ we have
\begin{equation*}
  \begin{split}
    \exp
    \left[
      \int_{z_0}^t p(u)\,dt
    \right]
    &= \exp
    \left[
      \int_{z_0}^t \left(-\frac1u\right)\,du + \int_{z_0}^t \left[ p(u) +
      \frac1u \right]\, du
    \right] \\
    &= \exp
    \left[
      -\ln t + \ln z_0 + \int_{z_0}^t \left[ p(u) + \frac1u \right]\,
      du
    \right] \\
    &\sim z_0 e^a t^{-1}, \text{ where } a := \int_{z_0}^0 \left[ p(u)
    + \frac1u \right]\, du,
  \end{split}
\end{equation*}
whereas, using~\eqref{eq:5} and~\eqref{eq:3},
\begin{equation*}
  g^{[s]}(t) \sim \frac{r^{[s]}(t)}t \sim T_2 r_2^{[s]} t;
\end{equation*}
thus the integrand in~\eqref{eq:6} has no singularity at~$t=0$.

Now we obtain the particular solution of interest, using the boundary
condition $\mu^{[s]}(z) \sim t_1^s T_1 z$ as $z \to 0$.  We find the
constant~$\beta_s$ is specified in terms of~$z_0$ as
\begin{equation}
   \label{eq:beta_s}
   \beta_s = z_0 e^a t_1^s T_1.
\end{equation}
\begin{remark}
  \label{rem:two-sided-closed}
  One can check for each very simple family that
\begin{equation*}
 \Phi'(t) = \frac{a_{0} + a_{1}}{1+a_{0}t} \Phi(t),
\end{equation*}
and for any simply generated family that
\begin{equation*}
  T'(z) = \frac{\Phi(T(z))}{1-\frac{\Phi'(T(z))
      T(z)}{\Phi(T(z))}}.
\end{equation*}
Thus
\begin{equation*}
  p(z) = - \frac{\Phi(T(z))}{T(z)
    [1-a_{1}T(z)]} \left[1+2a_{0} T(z) +
    \frac{a_{1} T(z)}
    {1-\frac{(a_{0}+a_{1})T(z)}{1+a_{0}T(z)}}\right].
\end{equation*}
This leads to
\begin{align*}
  \int_{z_0}^{z} p(t)\,dt &= - \int_{T(z_0)}^{T(z)} \frac1{T
    (1-a_{1}T)}{\left[1+2a_{0} T + 
    \frac{a_{1} T} {1-\frac{(a_{0}+a_{1})T}{1+a_{0}T}}\right]
    \left[1-\frac{(a_{0}+a_{1})T}{1+a_{0}T}\right]}\,dT \\
  & = -\int_{T(z_0)}^{T(z)} \Big(\frac{1}{T} + \frac{a_{0}}{1+a_{0}T}
  + \frac{a_{1}}{1-a_{1} T}\Big)\,dT \\
  &= \ln\left[\frac{1-a_{1}
      T(z)}{T(z) (1+a_{0}T(z))}\right] - \ln
  \left[
    \frac{1-a_1T(z_0)}{T(z_0)(1+a_0T(z_0))}
  \right]
\end{align*}
and finally,
again using the boundary conditions on $\mu^{[s]}(z)$ as $z \to 0$,
to the
following explicit form of~\eqref{eq:6.1}:
\begin{equation*}
  \muz{s} = \frac{T(z) [1+a_{0} T(z)]}{1-a_{1}T(z)} \left\{ \int_0^{z}
  g^{[s]}(t)  \frac{1-a_{1}T(t)}{T(t) [1+a_{0}T(t)]}\,dt + t_1^s \right\}.
\end{equation*}
\end{remark}

\subsection{One-sided destruction}
\label{sec:one-sided-destr-1}
Here, the cost of cutting down a very simple tree of size~$n$, 
call it $Y_n$, satisfies the distributional recurrence
\begin{equation}
  Y_n \eqinlaw Y_{K_n} + t_n, \quad n \geq 2; \qquad Y_1 = t_1,
  \label{eqno1}
\end{equation}
where $t_n$, for $n \geq 2$, is the toll for cutting an edge from a
tree of size~$n$ and the splitting probabilities are given by
$p_{n,k}$ at~\eqref{eq:p_nk}.

Defining $\mu_{n}^{[s]} := \E Y_{n}^{s}$, one obtains from equation
\eqref{eqno1} by conditioning on $K_{n}$ the recurrence relation
\begin{equation}
   \label{eqno2}
   \mu_{n}^{[s]} = \sum_{k=1}^{n-1} p_{n,k} \mu_{k}^{[s]} + r_{n}^{[s]}, \quad 
   n \ge 2, 
\end{equation}
where
\begin{equation*}
   r_{n}^{[s]} := \sum_{\substack{s_{1}+s_{2}=s, \\ s_{2} < s}} \binom{s}{s_{1}} t_{n}^{s_{1}}
   \sum_{k=1}^{n-1} p_{n,k} \mu_{k}^{[s_{2}]},
\end{equation*}
and $\mu_{1}^{[s]} = t_{1}^{s}$.
Using the same notation as in Section~\ref{sec:two-sided-destr}, we
obtain the following differential equation by
multiplying~\eqref{eqno2} by $(n-1)T_{n}z^{n}$ and summing over $n \ge
2$:
\begin{equation*}
   z \partial_{z} \mu^{[s]}(z) - \mu^{[s]}(z) = T(z) \big(a_{1} z \partial_{z} \mu^{[s]}(z)
   + a_{0} \mu^{[s]}(z)\big) + r^{[s]}(z),
\end{equation*}
where
\begin{equation}
  \label{eq:28}
   r^{[s]}(z) := \sum_{\substack{s_{1}+s_{2}=s, \\ s_{2} < s}} \binom{s}{s_{1}}
   t^{\odot s_{1}}(z) \odot \big[ T(z) \big(a_{1} z \partial_{z}
   \mu^{[s_{2}]}(z)
   + a_{0} \mu^{[s_{2}]}(z)\big)\big].
\end{equation}
This can be written as
\begin{equation}
   \partial_{z} \mu^{[s]}(z) + p(z) \mu^{[s]}(z) = g^{[s]}(z),
   \label{eqno3}
\end{equation}
with
\begin{equation}
  \label{eq:26}
   p(z) := - \frac{1+a_{0} T(z)}{z[1-a_{1} T(z)]} \quad \text{and} \quad
   g^{[s]}(z) := \frac{r^{[s]}(z)}{z [1-a_{1} T(z)]}.
\end{equation}

One can check that for each very simple family, $p(z) = - \partial_{z}
\ln(T(z))$, so that we obtain as general solution of the first order linear
differential equation \eqref{eqno3}:
\begin{equation*}
   \mu^{[s]}(z) = T(z) \int_{0}^{z} \frac{g^{[s]}(t)}{T(t)} dt + C \, T(z),
\end{equation*}
and finally by adapting to the initial condition $\left.\partial_{z}
  \mu^{[s]}(z)\right|_{z=0} = T_{1} \mu_{1}^{[s]} = T_{1} t_{1}^{s}$,
that the integration constant is given as $C=t_{1}^{s}$. Therefore, we
  get
\begin{equation}
   \label{eqno4}
   \mu^{[s]}(z) = T(z) \int_{0}^{z} \frac{g^{[s]}(t)}{T(t)} dt + t_{1}^{s} T(z).
\end{equation}

\section{Two-sided destruction}
\label{sec:toll-nalpha}

We begin by obtaining a singular expansion for $p(z)$ at~\eqref{eq:4}.
Using~\eqref{eq:9} and~\eqref{eq:10} in~\eqref{eq:4} we get
\begin{equation}
  \label{eq:38}
  p(z) \sim - \frac{\rho^{-1}}{2} Z^{-1} + Z^{-1/2}\A + \A.
\end{equation}
Integrating this singular expansion term-by-term~\cite[Theorem~7]{FFK},
\begin{equation*}
  \int_{z_0}^z p(t)\,dt \sim - \frac{1}{2} \ln Z^{-1} + \A + Z^{1/2}\A.
\end{equation*}
Thus
\begin{equation}
  \label{eq:11}
  \exp\left[ - \int_{z_0}^z p(t)\,dt \right] \sim \xi Z^{-1/2} + \A +
  Z^{1/2}\A,
\end{equation}
where
\begin{equation}
  \label{eq:12}
  \xi := (1 - \rho^{-1}z_0)^{1/2} \exp\left[ - \int_{z_0}^\rho \left[
  p(t) + \frac{\rho^{-1}}2(1-\rho^{-1}t)^{-1} \right]\,dt \right].
\end{equation}
Taking the reciprocal of~\eqref{eq:11} gives
\begin{equation}
  \label{eq:13}
  \exp \left[ \int_{z_0}^z p(t)\,dt \right] \sim \xi^{-1} Z^{1/2} + Z\A +
  Z^{3/2}\A.
\end{equation}

Let us now consider two-sided destruction with the toll $t_n =
n^\alpha$, with $\alpha > 0$. (Notice that the case $\alpha=0$ is
trivial since then the total cost of destruction is simply the number
of edges in the tree, which is always~$n-1$.)  The toll generating
function $t(z)$ is the generalized polylogarithm~$\Li_{-\alpha,0}(z)$,
which is amenable to singularity
analysis~\cite[Theorem~1]{MR2000a:05015}.

\subsection{Expectation}
\label{sec:mean}

Now we obtain a singular expansion for $r^{[1]}(z)$ defined
at~\eqref{eq:3}, recalling that $\mu^{[0]}(z) = T(z)$:
\begin{equation}
  \label{eq:31}
  r^{[1]}(z) = t(z) \odot [ a_1 zT'(z) T(z) + a_0 T^2(z) ].
\end{equation}
Using~\eqref{eq:7} we conclude that
\begin{equation*}
  T^2(z) \sim \A + Z^{1/2}\A,
\end{equation*}
and using~\eqref{eq:8} that
\begin{equation}
  \label{eq:13.1}
  a_1 z T'(z) T(z) + a_0 T^2(z) \sim \rho^{1/2} \frac{b}{2}
  Z^{-1/2}  + \A + Z^{1/2}\A.
\end{equation}
We will use the Zigzag algorithm of~\cite{FFK} to obtain a singular
expansion for $r^{[1]}(z)$.  We recall the use of the notation~$\N$ to
denote a generic power series in $1/n$, possibly different at each
occurrence.  By singularity analysis,
\begin{equation}
  \label{eq:33}
  \rho^n [z^n] [ a_1 z T'(z) T(z) + a_0 T^2(z) ] \sim \rho^{1/2}
  \frac{b}2 \frac{n^{-1/2}}{\sqrt\pi} + n^{-3/2} \N. 
\end{equation}
Thus
\begin{equation}
  \label{eq:23}
  \rho^n [z^n] r^{[1]}(z) \sim  \rho^{1/2} \frac{b}2
  \frac{n^{\alpha-\frac12}}{\sqrt\pi} +   n^{\alpha-\frac32} \N. 
\end{equation}

Until further notice, suppose $\alpha \not\in \{ \frac12, \frac32,
\ldots \}$.  Then a compatible singular
expansion for $r^{[1]}(z)$ at~\eqref{eq:31} is obtained as
\begin{equation}
  \label{eq:30}
  r^{[1]}(z) \sim \rho^{1/2}  b \frac{\Gamma(\alpha+\frac12)}{2
  \sqrt\pi} Z^{-\alpha-\frac12} + Z^{-\alpha+\frac12}\A + \A.
\end{equation}
Recalling~\eqref{eq:5} and~\eqref{eq:10} we have
\begin{equation}
  \label{eq:32}
  g^{[1]}(z) \sim \frac{\tau}{\rho}
  \frac{\Gamma(\alpha+\frac12)}{2\sqrt\pi} Z^{-\alpha-1} +
  Z^{-\alpha-\frac12}\A + Z^{-\alpha}\A + Z^{-1/2}\A + \A.
\end{equation}
Using this expansion and~\eqref{eq:13},
\begin{equation*}
  g^{[1]}(z) \exp\left[ \int_{z_0}^z p(t)\,dt \right] \sim
  \frac{\tau}{\xi \rho} 
  \frac{\Gamma(\alpha+\frac12)}{2\sqrt\pi} Z^{-\alpha-\frac12} +
  Z^{-\alpha}\A + Z^{-\alpha+\frac12}\A + \A + Z^{1/2}\A.
\end{equation*}
By~\eqref{eq:6} and Theorem~7 of~\cite{FFK}, we may integrate this
expansion term-by-term to get a complete singular 
expansion for~$A$.  If $\alpha \not\in \{1, 2, \dots\}$, we have
\begin{equation*}
  A^{[1]}(z) \sim \frac{\tau}{\xi} \frac{\Gamma(\alpha-\frac12)}{2\sqrt\pi}
  Z^{-\alpha + \frac12} + L_0 + Z^{-\alpha+1}\A + Z^{-\alpha +
  \frac32}\A + Z\A +  Z^{3/2}\A,
\end{equation*}
where $L_0$ is a constant.  [The value of~$L_0$ is immaterial unless
$0 < \alpha < 1/2$, in which case see~\eqref{eq:40}.] On the other
hand, if $\alpha \in \{1, 2, \dots \}$, a logarithmic term appears
upon integration, so that
\begin{equation*}
  A^{[1]}(z) \sim \frac{\tau}{\xi} \frac{\Gamma(\alpha-\frac12)}{2\sqrt\pi}
  Z^{-\alpha + \frac12} + Z^{-\alpha+1}\A + Z^{-\alpha +
  \frac32}\A +  K_0 \ln Z^{-1},
\end{equation*}
where $K_0$ is a constant.
Combining these expansions
with~\eqref{eq:11}, we finally obtain [recalling~\eqref{eq:6.1}]
\begin{equation}
  \label{eq:14}
  \muz{1} \sim \tau \frac{\Gamma(\alpha-\frac12)}{2\sqrt\pi}
  Z^{-\alpha} + L_0 \xi Z^{-1/2} + Z^{-\alpha+\frac12}\A +
  Z^{-\alpha+1}\A + \A + Z^{1/2}\A
\end{equation}
when $\alpha \not\in \{1,2,\dots\}$ and
\begin{equation}
  \label{eq:14a}
  \muz{1} \sim \tau \frac{\Gamma(\alpha-\frac12)}{2\sqrt\pi}
  Z^{-\alpha} + Z^{-\alpha+\frac12}\A +  Z^{-\alpha+1}\A +
  Z^{-1/2}(\ln Z^{-1})\A + (\ln Z^{-1})\A
\end{equation}
when $\alpha \in \{1,2, \dots\}$.  Note that the remainder
in~\eqref{eq:14a} is $O(|Z|^{-\alpha + \frac12})$ unless $\alpha=1$,
in which case it is $O(|Z|^{-\frac12} \ln Z^{-1}) = O(|Z|^{-\frac12 -
  \epsilon})$ for any $\epsilon > 0$.

When $\alpha > 1/2$ and $\alpha \not\in \{1,2,\dots\}$, by
singularity analysis we have
\begin{equation*}
  \rho^n \mun{1} T_n \sim \tau
  \frac{\Gamma(\alpha-\frac12)}{2\sqrt\pi\Gamma(\alpha)} n^{\alpha-1}
  + n^{\alpha-\frac32}\N + n^{\alpha-2}\N + n^{-1/2}\N,
\end{equation*}
so that, recalling~\eqref{eq:15.1} and~\eqref{eq:15},
\begin{equation*}
  \mun{1} \sim \sigma
  \frac{\Gamma(\alpha-\frac12)}{\sqrt2\Gamma(\alpha)}
  n^{\alpha+\frac12} + n^\alpha\N + n^{\alpha-\frac12}\N + n\N.
\end{equation*}
When $\alpha \in \{1, 2, \dots \}$, starting from~\eqref{eq:14a} and
the note following that display, we can similarly derive the expansion
\begin{equation*}
  \mun{1} = \sigma \frac{\Gamma(\alpha-\frac12)}{\sqrt2\Gamma(\alpha)}
  n^{\alpha+\frac12} + O(n^{\alpha}) + O(n \log n).
\end{equation*}

When $0 < \alpha < 1/2$, a similar computation yields
\begin{equation}
  \label{eq:16}
  \mun{1} \sim \frac{L_0 \xi}{c\sqrt\pi} n + \sigma
  \frac{\Gamma(\alpha-\frac12)}{\sqrt2\Gamma(\alpha)} n^{\alpha+\frac12} +
  n^\alpha\N + n^{\alpha-\frac12}\N + \N,
\end{equation}
where
\begin{equation}
  \label{eq:40}
  L_0 := \int_0^\rho g^{[1]}(t) \exp\left[ \int_{z_0}^t p(u)\,du
  \right] \,dt  + \beta_1,
\end{equation}
with $p$ and $g^{[1]}$ defined at~\eqref{eq:4} and~\eqref{eq:5},
respectively, and $\xi$ and~$\beta_1$ at~\eqref{eq:12}
and~\eqref{eq:beta_s}, respectively.

When $\alpha \in \{\frac32, \frac52, \dots \}$, one can check that
logarithmic terms appear in the singular expansion compatible
with~\eqref{eq:23} but the lead-order term and asymptotic order of the
remainder are unchanged. Indeed, now
\begin{equation}
  \label{eq:24}
  \muz{1} = \tau \frac{\Gamma(\alpha-\frac12)}{2\sqrt\pi}
  Z^{-\alpha} +   O(|Z|^{-\alpha + \frac12})
\end{equation}
and consequently
\begin{equation*}
  \mun{1} = \sigma\frac{\Gamma(\alpha-\frac12)}{\sqrt2\Gamma(\alpha)}
  n^{\alpha+\frac12} + O(n^{\alpha}).
\end{equation*}

Finally we consider $\alpha=1/2$.  Now, a compatible singular
expansion for~\eqref{eq:23} is
\begin{equation*}
  r^{[1]}(z) \sim \frac{\rho^{1/2}b}{2\sqrt\pi} Z^{-1} + (\log Z)\A +
  \A.
\end{equation*}
Proceeding as in the case $\alpha \ne 1/2$ we have the singular
expansions 
\begin{equation*}
  \begin{split}
    g^{[1]}(z) \sim \frac{\tau}{2 \rho \sqrt\pi} Z^{-3/2} + Z^{-1}\A +
    (Z^{-1/2} \log Z)\A + Z^{-1/2}\A + (\log Z)\A, \\
    g^{[1]}(z) \exp\left[ \int_{z_0}^z p(t)\,dt \right] \sim
    \frac{\tau}{\xi\rho}
    \frac1{2\sqrt\pi} Z^{-1} + Z^{-\frac12}\A + (\log Z)\A + \A +
    (Z^{1/2}\log Z)\A, \\
    A^{[1]}(z) \sim \frac{\tau}{\xi} \frac{1}{2\sqrt\pi} \ln Z^{-1} + L_1 +
    Z^{1/2}\A + Z\A + (Z \log Z)\A + (Z^{3/2}\log Z)\A,
  \end{split}
\end{equation*}
where
\begin{equation}
  \label{eq:39}
  L_1 := \int_0^\rho \left\{ g^{[1]}(t) \exp\left[ \int_{z_0}^t p(u)\,du
  \right] - \frac{\tau}{\xi\rho} \frac1{2\sqrt\pi} (1-\rho^{-1}t)
  \right\} \,dt + \beta_1.
\end{equation}
This leads to
\begin{equation*}
  \muz{1} \sim \frac{\tau}{2\sqrt\pi} Z^{-1/2} \ln Z^{-1} + \xi L_1
  Z^{-1/2} + (\log Z)\A + \A + (Z^{1/2} \log Z)\A + Z^{1/2}\A,
\end{equation*}
so that by singularity analysis and~\eqref{eq:15.1} we have
\begin{equation}
  \label{eq:35}
  \begin{split}
\mun{1} \sim \frac{\sigma}{\sqrt{2\pi}} n \ln n + \left[ \frac{\xi
      L_1}{c \sqrt\pi} + \frac{\sigma}{\sqrt{2\pi}}(\gamma + 2\ln 2) \right]
  n \\
  {} + (n^{1/2} \log n) \N + (\log n)\N + \N,
\end{split}
\end{equation}
where $\sigma$ is defined at~\eqref{eq:15}.

\subsection{Higher moments and limiting distributions}
\label{sec:higher-moments}

We proceed to higher moments.  We will consider separately the cases
$\alpha > 1/2$, $\alpha < 1/2$, and $\alpha = 1/2$.  We present the
details for $\alpha > 1/2$ and sketch the main ideas for the other
cases. Throughout $\alpha' := \alpha + \frac12$.

\begin{proposition}
  \label{prop:gt12}
  Let $\alpha > 1/2$ and $\epsilon > 0$.  Then
  \begin{equation*}
    \muz{s} = C_s Z^{-s\alpha' + \frac12} + O(|Z|^{-s\alpha' +
      \frac12 + q}),
  \end{equation*}
  where
  \begin{equation*}
    q := 
    \begin{cases}
      \min\{\alpha-\frac12,\frac12\} & \text{\textup{if} $\alpha \ne 1$} \\
      \frac12 - \epsilon & \text{\textup{if} $\alpha = 1$}
    \end{cases}
  \end{equation*}
  with
  \begin{equation*}
    C_1 = \tau \frac{\Gamma(\alpha-\frac12)}{2\sqrt\pi},
  \end{equation*}
  and, for $s \geq 2$,
  \begin{equation}
    \label{eq:17}
    C_s = \rho^{-1/2}b^{-1} \left[ \frac{1}{s\alpha'-1}
      \sum_{k=1}^{s-1} \binom{s}{k}
      \left(k\alpha'-\frac12 \right) C_{k} C_{s-k} + s \tau
    \frac{\Gamma(s\alpha'-1)}{\Gamma((s-1)\alpha'-\frac12)}
      C_{s-1}\right].
  \end{equation}
\end{proposition}
\begin{proof}
  The proof is by induction on $s$.  The claim is true for $s=1$
  by~\eqref{eq:14}, \eqref{eq:14a}, and~\eqref{eq:24}.  Suppose $s \geq
  2$.  We analyze each term in the sum for $r^{[s]}(z)$
  at~\eqref{eq:3}.

  If both $s_2$ and $s_3$ are nonzero, then by the induction hypothesis,
  \begin{equation*}
      z\partial_z \muz{s_2} = C_{s_2}\left(s_2\alpha' -
      \frac12\right)
      Z^{-s_2\alpha' - \frac12} + O(|Z|^{-s_2\alpha' - \frac12 + q}),
  \end{equation*}
  so that
  \begin{equation*}
    z\left( \partial_z \muz{s_2} \right) \muz{s_3} = C_{s_2} C_{s_3}
      \left(s_2\alpha' - \frac12\right) Z^{-(s_2+s_3)\alpha'} +
      O(|Z|^{-(s_2 + s_3)\alpha' + q}).
  \end{equation*}
  Also, $\muz{s_2} \muz{s_3} = O(|Z|^{-(s_2+s_3)\alpha'+1}).$  Hence
  \begin{multline*}
    a_1 z \left(\partial_z \muz{s_2} \right) \muz{s_3} + a_0 \muz{s_2}
    \muz{s_3} \\
    = a_1 C_{s_2} C_{s_3} \left(s_2\alpha' - \frac12 \right)
    Z^{-(s_2+s_3)\alpha'} + O(|Z|^{-(s_2+s_3)\alpha' + q}).
  \end{multline*}
  Taking the Hadamard product of this expansion with $t^{\odot
  s_1}(z)$ (using the Zigzag algorithm again)
  gives the contribution of such terms to $r^{[s]}(z)$ as
  \begin{equation*}
    \binom{s}{s_1,s_2,s_3} a_1 C_{s_2} C_{s_3} \left(s_2 \alpha' -
    \frac12 \right) 
    \frac{\Gamma(s\alpha' - \frac{s_1}2)}{\Gamma((s_2+s_3)\alpha')}
    Z^{-(s\alpha' - \frac{s_1}{2})} + O(|Z|^{-(s\alpha' -
    \frac{s_1}{2}) + q}).
  \end{equation*}
  Notice that if $s_1 \ne 0$ the contribution is
  $O(|Z|^{-s\alpha'+\frac12})$.

  Next consider the case when $s_2$ is nonzero but $s_3 = 0$.  By the
  induction hypothesis and the singular expansion of $T$
  at~\eqref{eq:7},
  \begin{equation*}
      z\left( \partial_z \muz{s_2} \right) T(z) = \tau C_{s_2}
      \left(s_2\alpha'-\frac12\right) Z^{-s_2\alpha'-\frac12} +
      O(|Z|^{-s_2\alpha'-\frac12+q}).
  \end{equation*}
  Also $\muz{s_2}T(z) = O(|Z|^{-s_2\alpha' + \frac12})$.  Hence
  \begin{multline*}
    a_1z\left( \partial_z \muz{s_2} \right) T(z) + a_0 \muz{s_2}T(z)\\
    =  C_{s_2} \left(s_2\alpha'-\frac12\right)
    Z^{-s_2\alpha'-\frac12} + O(|Z|^{-s_2\alpha'-\frac12+q}).
  \end{multline*}
  Taking the Hadamard product of this singular expansion with
  $t^{\odot s_1}(z)$ we get that the contribution to $r^{[s]}(z)$ from
  such terms is
  \begin{equation*}
    \binom{s}{s_1}  C_{s_2} \left(s_2
    \alpha'-\frac12\right)
    \frac{\Gamma(s\alpha'-\frac{s_1}2+\frac12)}{\Gamma(s_2\alpha'+\frac12)}
    Z^{-(s\alpha'-\frac{s_1}2+\frac12)} +
    O(|Z|^{-(s\alpha'-\frac{s_1}2+\frac12)+q}).
  \end{equation*}
  Notice that $s_1 \geq 1$ and that when $s_1 > 1$ the contribution of
  such terms is $O(|Z|^{-s\alpha'+\frac12})$.

  We move on to the case when $s_2=0$ but $s_3$ is
  nonzero. By the
  induction hypothesis,~\eqref{eq:7}, and~\eqref{eq:8}, we have $T(z)
  \muz{s_3} =
  O(|Z|^{-s_3\alpha'+\frac12})$ and $zT'(z) \muz{s_3} =
  O(|Z|^{-s_3\alpha'})$.  Thus
  \begin{equation*}
    a_1 zT'(z) \muz{s_3} + a_0 T(z) \muz{s_3} = O(|Z|^{-s_3\alpha'}).
  \end{equation*}
  Taking the Hadamard product with $t^{\odot s_1}(z)$ we see
  (recalling $s_1 \geq 1$) that the contribution to $r^{[s]}(z)$ from
  these terms is $O(|Z|^{-s\alpha'+\frac12})$.

  Finally we consider the case when $s_2=s_3=0$.  In this case,
  using~\eqref{eq:13.1} it is easy to verify that the contribution to
  $r^{[s]}(z)$ from this term is $O(|Z|^{-s\alpha'+\frac12})$.

  Summing all the contributions we see that
  \begin{equation}
    \label{eq:37}
    r^{[s]}(z) = D_s Z^{-s\alpha'} + O(|Z|^{-s\alpha'+q}),
  \end{equation}
  where
  \begin{equation*}
    D_s := a_1 \left[ \sum_{k=1}^{s-1} \binom{s}{k}
    \left(k\alpha'-\frac12 \right) C_{k} C_{s-k} + s \tau 
    \frac{[(s-1)\alpha' -
    \frac12]\Gamma(s\alpha')}{\Gamma((s-1)\alpha'+\frac12)} C_{s-1}
    \right]
  \end{equation*}
  Thus, using~\eqref{eq:5} and~\eqref{eq:10},
  \begin{equation*}
    g^{[s]}(z) = \rho^{-3/2} a_1^{-1} b^{-1} D_s Z^{-s\alpha'-\frac12} +
    O(|Z|^{-s\alpha'-\frac12+q}),
  \end{equation*}
  whence, using~\eqref{eq:13},
  \begin{equation*}
    g^{[s]}(z) \exp\left[ \int_{z_0}^z p(t)\,dt \right] = \xi^{-1}
    \rho^{-3/2} a_1^{-1} b^{-1} D_s Z^{-s\alpha'} + O(|Z|^{-s\alpha'+q}).
  \end{equation*}
  To get $A^{[s]}(z)$ at~\eqref{eq:6} we integrate this singular expansion,
  noting that since $s \geq 2$ and $\alpha' > 1$, we have $s\alpha'
  > 2$.  Hence
  \begin{equation*}
    A^{[s]}(z) = \xi^{-1} \rho^{-1/2} a_1^{-1} b^{-1} \frac{D_s}{s\alpha'-1}
    Z^{-s\alpha'+1} + O(|Z|^{-s\alpha'+1+q}).
  \end{equation*}
  Now by~\eqref{eq:6.1} and~\eqref{eq:11},
  \begin{equation*}
    \muz{s} = \rho^{-1/2} a_1^{-1} b^{-1} \frac{D_s}{s\alpha'-1}
    Z^{-s\alpha'+\frac12} + O(|Z|^{-s\alpha'+\frac12+q}).
  \end{equation*}
  Taking
  \begin{align*}
    C_s &= \rho^{-1/2} a_1^{-1} b^{-1} \frac{D_s}{s\alpha'-1} \\
    &= \frac{\rho^{-1/2}b^{-1}}{s\alpha'-1} \left[ \sum_{k=1}^{s-1}
    \binom{s}{k}
      \left(k\alpha'-\frac12 \right) C_{k} C_{s-k} + s \tau \frac{[
        (s-1)\alpha' -
        \frac12]\Gamma(s\alpha')}{\Gamma((s-1)\alpha'+\frac12)}
      C_{s-1}
    \right] \\
    &= \rho^{-1/2}b^{-1} \left[ \frac{1}{s\alpha'-1} \sum_{k=1}^{s-1}
    \binom{s}{k}
      \left(k\alpha'-\frac12 \right) C_{k} C_{s-k} + s \tau
    \frac{\Gamma(s\alpha'-1)}{\Gamma((s-1)\alpha'-\frac12)} C_{s-1}\right]
  \end{align*}
  completes the proof.
\end{proof}
Using singularity analysis we can now derive asymptotics for the
moments~$\mun{s}$.
\begin{theorem}
  \label{thm:gt12}
  Let $\alpha > 1/2$. Then, as $n \to \infty$,
  \begin{equation*}
    \sigma^{-s}n^{-s\alpha'}\mun{s} \to m_s,
  \end{equation*}
  where $\sigma^2 := \tau^2 \frac{\Phi''(\tau)}{\Phi(\tau)}$ and $m_s$
  (which does not depend on the very simple family) is given by
  \begin{equation*}
    m_1 = \frac{\Gamma(\alpha-\frac12)}{\sqrt2 \Gamma(\alpha)}
  \end{equation*}
  and, for $s \geq 2$,
  \begin{equation}
    \label{eq:18}
    m_s = \frac{1}{4\sqrt\pi}
    \sum_{k=1}^{s-1}
    \binom{s}{k} \frac{\Gamma(k\alpha'-\tfrac12)
    \Gamma((s-k)\alpha'-\tfrac12)}{\Gamma(s\alpha'-\frac12)} 
    m_k m_{s-k}  +
    \frac{s\Gamma(s\alpha'-1)}{\sqrt2\Gamma(s\alpha'-\frac12)}
    m_{s-1}.
  \end{equation}
\end{theorem}
\begin{proof}
  Using singularity analysis and Proposition~\ref{prop:gt12},
  \begin{equation*}
    \rho^n \mun{s} T_n = C_s
    \frac{n^{s\alpha'-\frac32}}{\Gamma(s\alpha'-\frac12)} +
    O(n^{s\alpha'-\frac32-q}),
  \end{equation*}
  and using the asymptotics of $T_n$ at~\eqref{eq:15.1},
  \begin{equation*}
    \mun{s} = \frac{C_s}{c \Gamma(s\alpha'-\frac12)} n^{s\alpha'} +
    O(n^{s\alpha'-q}).
  \end{equation*}
  Then
  \begin{equation*}
    \sigma^{-s} n^{-s\alpha'} \mun{s} \to m_s,
  \end{equation*}
  where
  \begin{equation}
    \label{eq:19}
    m_s := \sigma^{-s} \frac{C_s}{c \Gamma(s\alpha'-\frac12)}.
  \end{equation}
  Thus, using $2\sqrt\pi c \sigma = \sqrt2\tau$,
  \begin{equation*}
    m_1 = \frac{C_1}{c\sigma\Gamma(\alpha)} =
    \frac{\Gamma(\alpha-\frac12)}{\sqrt2 \Gamma(\alpha)}
  \end{equation*}
  Using~\eqref{eq:17},~\eqref{eq:19}, and the identities
  \begin{equation*}
      c\rho^{-1/2}b^{-1} = \frac{1}{2\sqrt{\pi}}, \qquad
      \sigma^{-1} \tau \rho^{-1/2} b^{-1} = \frac1{\sqrt2}, \qquad
      \Gamma(x+1) = x\Gamma(x),
  \end{equation*}
  we obtain the following recurrence for $m_s$:
    \begin{equation}
      \label{eq:20}
    m_s = \frac{1}{2\sqrt\pi}
    \sum_{k=1}^{s-1}
    \binom{s}{k} \frac{\Gamma(k\alpha'+\tfrac12)
    \Gamma((s-k)\alpha'-\tfrac12)}{(s\alpha'-1)\Gamma(s\alpha'-\frac12)}
    m_k m_{s-k} \\ +
    \frac{s\Gamma(s\alpha'-1)}{\sqrt2\Gamma(s\alpha'-\frac12)} m_{s-1}.
  \end{equation}
  To obtain the form of the recurrence in~\eqref{eq:18},
  write~\eqref{eq:20} in the form
  \begin{equation*}
    m_s = \frac{1}{2\sqrt\pi} \sum_{k=1}^{s-1} e_{s,k}
    + \tilde{e}_s = \frac{1}{4\sqrt\pi} \sum_{k=1}^{s-1} (e_{s,k} + e_{s,s-k})
    + \tilde{e}_s
  \end{equation*}
  and simplify.
\end{proof}
\begin{remark}
  \label{rem:symmetry}
  In going from~\eqref{eq:20} to~\eqref{eq:18} we symmetrized by
  collecting coefficients of $m_k m_{s-k}$.  We might also have
  symmetrized from the start by choosing the splitting probabilities
  as
  \begin{equation*}
    \tilde{p}_{n,k} := \frac12( p_{n,k} + p_{n,n-k} ).
  \end{equation*}
  In the particular case of Cayley trees this leads to the same
  splitting probabilities as for the Union--Find recurrence studied
  in~\cite{MR81a:68049,MR89i:05024,FFK}.
\end{remark}

We can now show convergence in distribution via the method of moments.
\begin{theorem}
  \label{thm:dist12}
  Let $\alpha > 1/2$.
  Define $\sigma^2 := \tau^2\Phi''(\tau)/\Phi(\tau)$ and $\alpha' :=
  \alpha + \frac12$.  Then, as $n \to \infty$, 
  \begin{equation*}
    \sigma^{-1} n^{-\alpha'} X_n \toinlaw X^{(\alpha)},
  \end{equation*}
  with convergence of all moments, where $X^{(\alpha)}$ has the unique
  distribution whose $s$th moment $m_s \equiv m_s(\alpha)$ is given by
  \begin{equation*}
    m_1 = \frac{\Gamma(\alpha-\frac12)}{\sqrt2 \Gamma(\alpha)}
  \end{equation*}
  and for $s \geq 2$ by the recurrence~\eqref{eq:18}.
\end{theorem}
\begin{proof}
  One need only check that the $m_k$'s satisfy Carleman's condition.
  This has already been established in~\cite{fill03:_limit_catal}.
\end{proof}
\begin{remark}
  \label{rem:curious}
It is curious that $\sigma^{-1} n^{-\alpha'} X_n$ has the same limiting
distribution as 
\begin{equation*}
\sigma n^{-\alpha'} \sum_{v \in T} |T_v|^\alpha.  
\end{equation*}
Here $T$ is a random simply generated tree and $|T_v|$ denotes the
size of the tree rooted at a node~$v$.  This was established
in~\cite{sgtechreport}.
\end{remark}

For the case $0 < \alpha < 1/2$ it is convenient instead to consider the
random variable
\begin{equation*}
  \widetilde{X}_n := X_n - \mu n,\qquad \mu :=
  \frac{L_0\xi}{c\sqrt\pi}.
\end{equation*}
[Note that, by~\eqref{eq:16}, $\mu n$ is the lead term in the
asymptotics of $\E X_n$ when $\alpha < 1/2$.]  Using~\eqref{eq:1},
\begin{equation}
  \label{eq:21}
  \widetilde{X}_n \eqinlaw \widetilde{X}_{K_n} + \widetilde{X}_{n-K_n}^*
  + t_n, \quad n \geq 2; \qquad \widetilde{X}_1 = t_1 - \mu.
\end{equation}
Define $\tilde{\mu}_n^{[s]} := \E{\widetilde{X}_n^s}$ and
\begin{equation*}
  \tilde{\mu}^{[s]}(z) := \sum_{n \geq 1} \tilde{\mu}_n^{[s]} T_n z^n.
\end{equation*}
Then, in analogous fashion,~\eqref{eq:3}--\eqref{eq:6} hold with
$\muz{s}$ replaced by $\tilde{\mu}^{[s]}(z)$.

Observe that, by~\eqref{eq:14},
\begin{equation}
  \label{eq:22}
  \tilde{\mu}^{[1]}(z) \sim \tau \frac{\Gamma(\alpha-\frac12)}{2\sqrt\pi}
  Z^{-\alpha} + (Z^{-\alpha+\frac12} +
  Z^{-\alpha+1} + 1 + Z^{1/2})\A.
\end{equation}
We can use~\eqref{eq:22} and~\eqref{eq:3}--\eqref{eq:6} to show that
Proposition~\ref{prop:gt12} holds for $\alpha < 1/2$ with $\muz{s}$
replaced by $\tilde{\mu}^{[s]}(z)$ and $q$ changed to
$2\alpha-\epsilon$, for sufficiently small $\epsilon>0$.  It follows
then that $X_n - \mu n$ has (after scaling) a limiting distribution.
\begin{theorem}
  \label{thm:lt12}
  Let $\alpha < 1/2$.
    Define $\sigma^2 := \tau^2\Phi''(\tau)/\Phi(\tau)$ and $\alpha' :=
  \alpha + \frac12$.  Then, as $n \to \infty$,
  \begin{equation*}
    \sigma^{-1} n^{-\alpha'} (X_n - \mu n) \toinlaw X^{(\alpha)},
  \end{equation*}
  with convergence of all moments, where $X^{(\alpha)}$ has the unique
  distribution whose $s$th moment $m_s \equiv m_s(\alpha)$ is given
  for $s=1$ by
  \begin{equation*}
    m_1 = \frac{\Gamma(\alpha-\frac12)}{\sqrt2 \Gamma(\alpha)}
  \end{equation*}
  and for $s \geq 2$ by the recurrence~\eqref{eq:18}.
\end{theorem}

Finally we turn our attention to the case $\alpha=1/2$.  Now, we
define
\begin{equation*}
  \wX_n := X_n - \frac{\sigma}{\sqrt{2\pi}} n \ln n - \delta n \text{
  with }
  \delta := \frac{\xi L_1}{c \sqrt\pi} + \frac{\sigma}{\sqrt{2\pi}}
  (\gamma + 2 \ln 2),
\end{equation*}
with $L_1$ defined at~\eqref{eq:39}.
Then [cf.~\eqref{eq:1}]
\begin{equation*}
  \wX_n \eqinlaw \wX_{K_n} + \wX_{n-K_n}^* + t_{n,K_n}, \quad n \geq 2,
\end{equation*}
with $\wX_1 = 1 - \delta$ and
\begin{equation*}
  t_{n,k} := \frac{\sigma}{\sqrt{2\pi}}\left[ k \ln k + (n-k) \ln{(n-k)} -
  n \ln n + \frac{\sqrt{2\pi}}\sigma n^{1/2}\right].
\end{equation*}
As in the case $\alpha < 1/2$, it is easily checked
that~\eqref{eq:3}--\eqref{eq:6} hold with $\muz{s}$ replaced by
$\tilde{\mu}^{[s]}(z)$ and $r_n^{[s]}$ at~\eqref{eq:29} replaced by
\begin{equation}
  \label{eq:36}
  \tilde{r}_n^{[s]} :=
  \sums \binom{s}{s_1,s_2,s_3}  \sum_{k=1}^{n-1}
  p_{n,k} t_{n,k}^{s_1} \tilde{\mu}_{k}^{[s_2]}
  \tilde{\mu}_{n-k}^{[s_3]}.
\end{equation}
The limiting distribution is given by the following result.
\begin{theorem}
  \label{thm:moments-half}
  As $n \to \infty$,
  \begin{equation*}
    \sigma^{-s} n^{-s} \tilde{\mu}_n^{[s]} \to m_s,
  \end{equation*}
  where $m_0 = 1$, $m_1 = 0$, and for $s \geq 2$,
  \begin{equation*}
    m_s = \frac{\Gamma(s-1)}{2\sqrt\pi\Gamma(s-\frac12)} \sums
    \binom{s}{s_1, s_2, s_3} \left( \frac{1}{\sqrt{2\pi}}
    \right)^{s_1} m_{s_2} m_{s_3} J_{s_1, s_2, s_3},
  \end{equation*}
  with
  \begin{equation*}
    J_{s_1, s_2, s_3} := \int_0^1 [ x \ln x + (1-x) \ln(1-x)]^{s_1}
    x^{s_2 - \frac12} (1-x)^{s_3 - \frac32} \, dx.
  \end{equation*}
  Consequently
  \begin{equation*}
    \sigma^{-1} n^{-1} \wX_{n} \toinlaw X^{(1/2)},
  \end{equation*}
  where $X^{(1/2)}$ has the unique distribution whose $s$th moment is
  given by $m_s$.
\end{theorem}
\begin{proof}[Proof sketch]
  We provide an outline of the proof, leaving the details to the
  reader.  We claim that it is sufficient to show that
  \begin{equation}
    \label{eq:34}
    \rho^n \tilde{\mu}_n^{[s]} T_n = [C_s + o(1)]n^{s-\frac32},
  \end{equation}
  with $C_0 = c$, $C_1 = 0$, and for $s \geq 2$,
  \begin{equation*}
    C_s = \frac{1}{b \sqrt\rho}\frac{\Gamma(s-1)}{\Gamma(s-\frac12)}
    \sums \binom{s}{s_1, s_2, s_3} \left(\frac{\sigma}{\sqrt{2\pi}}
    \right)^{s_1} C_{s_2} C_{s_3} J_{s_1, s_2, s_3}.
  \end{equation*}
  Indeed, defining $m_s := \sigma^{-s} c^{-1} C_s$ and proceeding as
  in Theorem~\ref{thm:gt12} yields the claim.

  To show~\eqref{eq:34}, we proceed by induction.  The case~$s=0$ is
  easily checked, and the case $s=1$ follows from~\eqref{eq:35}.  For $s \geq
  2$ we use the induction hypothesis and approximation of sums by
  Riemann integrals in~\eqref{eq:36} to get
  \begin{equation*}
    \rho^n (n-1) T_n \tilde{r}_n^{[s]} \sim D_s n^{s-1},
  \end{equation*}
  where
  \begin{equation*}
    D_s := a_1 \sums \binom{s}{s_1, s_2, s_3}
    \left(\frac{\sigma}{\sqrt{2\pi}} \right)^{s_1} C_{s_2} C_{s_3}
    J_{s_1,  s_2, s_3}.
  \end{equation*}
  Since we know a priori that $\tilde{r}^{[s]}(z)$ is amenable to singularity
  analysis it follows that [cf.~\eqref{eq:37}]
  \begin{equation*}
    \tilde{r}^{[s]}(z) \sim \Gamma(s){D_s} Z^{-s}
  \end{equation*}
  and completing the computations as in the proof of
  Proposition~\ref{prop:gt12} yields the proof of~\eqref{eq:34}.
\end{proof}

\section{One-sided destruction}
\label{sec:one-sided-destr}

\subsection{Expectation\label{seco1}}

We study equation \eqref{eqno4} for the toll $t_{n} = n^{\alpha}$ with
$\alpha \ge 0$ and start by establishing a singular expansion for the
expectation $\mu^{[1]}(z)$.  
Since $\mu^{[0]}(z) = T(z)$, we have from~\eqref{eq:28} that
\begin{equation*}
   r^{[1]}(z) = t(z) \odot \big[a_{1} z T'(z) T(z) + a_{0} T^{2}(z)\big],
\end{equation*}
which has already been considered in
Section~\ref{sec:mean}.
In the remaining part of Section~\ref{seco1}, we suppose now 
$\alpha \not\in \{ \frac12, \frac32, \ldots \} \cup \{0,1,2,\ldots
\}$. (The complementary cases are covered in the proof of
Theorem~\ref{theo1}.)  Then a compatible singular expansion for
$r^{[1]}(z)$ is available at~\eqref{eq:30}.
This leads to the expansion~\eqref{eq:32} for $g^{[1]}(z)$ 
and consequently, using~\eqref{eq:27}, to 
\begin{equation*}
   \frac{g^{[1]}(z)}{T(z)} \sim \frac{\Gamma(\alpha+\frac{1}{2})}{2 \rho \sqrt{\pi}}
   Z^{-\alpha-1} + Z^{-\alpha-\frac{1}{2}} \A
   + Z^{-\alpha} \A + Z^{-{1}/{2}} \A + \A.
\end{equation*}

Integrating the last expression gives the singular expansion 
\begin{equation*}
   \int_{0}^{z} \frac{g^{[1]}(t)}{T(t)} dt \sim 
   \frac{\Gamma(\alpha+\frac{1}{2})}{2 \alpha \sqrt{\pi}}
   Z^{-\alpha} + Z^{-\alpha+\frac{1}{2}} \A
   + Z^{-\alpha+1} \A + \A + Z^{{1}/{2}} \A.
\end{equation*}
Now using~\eqref{eqno4}, we obtain easily the desired 
expansion
for $\mu^{[1]}(z)$:
\begin{align}
   \mu^{[1]}(z) & = T(z) \int_{0}^{z} \frac{g^{[1]}(t)}{T(t)} dt +
   t_{1} T(z) \notag \\ 
   & \sim \frac{\tau \Gamma(\alpha+\frac{1}{2})}{2 \alpha \sqrt{\pi}}
   Z^{-\alpha} + Z^{-\alpha+\frac{1}{2}} \A + Z^{-\alpha+1} \A  + \A +
   Z^{{1}/{2}} \A. \label{eq:25}
\end{align}
Via singularity analysis, we thus get 
the following expansion
for the coefficients:  
\begin{equation*}
   \rho^n [z^{n}] \mu^{[1]}(z) = 
   \rho^{n} \mu_{n}^{[1]} T_{n} \sim \frac{\tau
   \Gamma(\alpha+\frac{1}{2})}{2 \sqrt{\pi} \, \Gamma(\alpha+1)}
   n^{\alpha-1} + n^{\alpha-\frac{3}{2}} \N + n^{\alpha-2} \N  +
   n^{-{3}/{2}} \N, 
\end{equation*}
which together with~\eqref{eq:15.1} 
yields 
the full asymptotic expansion
\begin{equation}
   \mu_{n}^{[1]} \sim \frac{\sigma
   \Gamma(\alpha+\frac{1}{2})}{\sqrt{2} \, \Gamma(\alpha+1)} 
   n^{\alpha+\frac{1}{2}} + n^{\alpha} \N + n^{\alpha-\frac{1}{2}} \N
   + \N,    \label{eqno5}
\end{equation}
with $\sigma$ defined at~\eqref{eq:15}.

\subsection{Higher moments and limiting distributions}
We state the main result of this section: 
\begin{theorem}
  \label{theo1}
  Let $\alpha \geq 0$. Define $\sigma := \tau
  \sqrt{\frac{\Phi''(\tau)}{\Phi(\tau)}}$ and $\alpha' := \alpha +
  \frac{1}{2}$. Then, for toll function $t_{n} = n^{\alpha}$, the
  moments 
  $\mu_n^{[s]} := \E Y_{n}^{s}$ satisfy the following
  asymptotic expansion as 
  $n \to \infty$: 
  \begin{equation*}
     \mu_{n}^{[s]} = \frac{s! \sigma^{s}}{2^{{s}/{2}}}
         \prod_{j=1}^{s}\frac{\Gamma(j \alpha')}{\Gamma(j\alpha'+\frac{1}{2})}
         n^{s \alpha'} + O\big(n^{s \alpha'-q }\big),
  \end{equation*}
  with
  \begin{equation*}
    q := 
    \begin{cases}
      \frac12 - \epsilon & \text{\textup{if} $\alpha \in \{0, 1/2\}$}
      \\
      \min\{\alpha, 1/2\} & \textup{otherwise},
    \end{cases}
  \end{equation*}
  where $\epsilon > 0$ is arbitrarily small.  Thus the normalized
  random 
  variable~$Y_{n}$ converges weakly to a random
  variable~$Y^{(\alpha)}$: 
  \begin{equation*}
     \sigma^{-1} n^{-\alpha'} Y_{n} \toinlaw Y^{(\alpha)},
  \end{equation*}
  where $Y^{(\alpha)}$ has the unique distribution with (for $s \geq
  1$) $s$th moment
  \begin{equation*}
     m_{s} = \frac{s!}{2^{{s}/{2}}}
         \prod_{j=1}^{s} \frac{\Gamma(j
         \alpha')}{\Gamma(j\alpha'+\frac{1}{2})}.
  \end{equation*}
  In particular when $\alpha = 0$ (i.e., $t_n \equiv 1$),
  $\sigma n^{-1/2} Y_n$ converges weakly to a standard
  Rayleigh distributed random variable~$Y^{(0)}$ with density 
  \begin{equation*}
    f(y) = y e^{-y^2/2}, \quad y \geq 0.
  \end{equation*}
  In this case the asymptotics of $\mun{s}$ can be sharpened to
  \begin{equation*}
     \mu_{n}^{[s]} = \frac{s! \sigma^s \sqrt\pi}{2^{s/2}
         \Gamma(\frac{s+1}{2})}
         n^{{s}/{2}} \, \left[1 + O\left(\frac{\log
         n}{\sqrt{n}}\right) \right]. 
  \end{equation*}
\end{theorem}
\begin{proof}
  We use induction on~$s$.  We begin with $\alpha > 0$.  
  Observe
  that it is sufficient to show that the generating
  functions~$\mu^{[s]}(z)$ admit the asymptotic
  expansions~\eqref{eqno6} around their dominant singularities at
  $z=\rho$. 
  Then, using singularity analysis, the claim
  follows. What we will show is that
\begin{equation}
   \label{eqno6}
   \mu^{[s]}(z) = 
          \frac{s! \sigma^{s-1} \tau}{2^{\frac{s+1}{2}} (s
      \alpha'-\frac{1}{2})  
      \sqrt{\pi}} \gamma_s Z^{-s \alpha'+\frac{1}{2}} +
      O\big(|Z|^{-s\alpha'+ \frac12 + q}\big),
\end{equation}
where
\begin{equation*}
  \gamma_s := \frac{\prod_{j=1}^{s} \Gamma(j \alpha')}{\prod_{j=1}^{s-1}
      \Gamma(j \alpha'+\frac{1}{2})}.
\end{equation*}

First 
we consider $s=1$, where we immediately obtain from the
full expansion~\eqref{eq:25} that~\eqref{eqno6} is true for all
$\alpha \not\in \{ \frac12, \frac32, \ldots \} \cup \{1,2,\ldots \}$.
If on the other hand $\alpha \in \{\frac{1}{2}, \frac32, \ldots \}
\cup \{1,2,\ldots \}$, then, repeating the computations of
Section~\ref{seco1}, it is easily seen 
that logarithmic terms
appear in the expansion of $\mu^{[1]}(z)$. But apart from the case
$\alpha=\frac{1}{2}$, they don't 
have an influence on the main
term or on the asymptotic growth order of the second-order 
term. If
$\alpha = \frac{1}{2}$, one observes that the general formula for the
main term holds, but the bound for the remainder term is different:\
$O(|\log Z^{-1}|)$, not $O(1)$.  Summarizing these cases, the
expansion~\eqref{eqno6}
holds for $s=1$. 

Next we assume that \eqref{eqno6} holds for all $1 \le s_{2} < s$ with
a given $s > 1$.  From~\eqref{eqno6} follows the expansion 
\begin{equation*}
   \partial_{z} \mu^{[s_{2}]}(z) = \frac{s_{2}! \sigma^{s_{2}-1} \tau}
   {2^{\frac{s_{2}+1}{2}} \sqrt{\pi} \rho}
   \gamma_{s_2} Z^{-s_{2} \alpha' - \frac{1}{2}}
   + O\big(|Z|^{-s_{2} \alpha' - \frac{1}{2} + q }\big),
\end{equation*}
which holds for all $1 \le s_{2} < s$.  Together with $\mu^{[0]}(z) =
T(z)$ and $a_1 \tau = 1$, this gives the following singular
expansion: 
\begin{multline*}
   T(z) \big(a_{1} z \partial_{z} \mu^{[s_{2}]}(z) + a_{0}
   \mu^{[s_{2}]}(z)\big) = \\ 
   \begin{cases}
   \frac12 b \rho^{1/2}  Z^{-{1}/{2}}
   + O(1), & s_{2}=0 \\
   \frac{s_{2}! \sigma^{s_{2}-1} \tau}{2^{(s_{2}+1)/{2}}
   \sqrt{\pi}} 
   \gamma_{s_2}
   Z^{-s_{2} \alpha' - \frac{1}{2}} + O\big(|Z|^{-s_{2} \alpha'
   - \frac12 + q}\big),
   & 1 \leq s_2 < s.
   \end{cases}
\end{multline*}

Under the assumptions $s_{1}+s_{2}=s$ and $s_{2}<s$, we get via
singularity analysis the expansion 
\begin{multline*}
   \rho^{n} [z^{n}] \binom{s}{s_{1}} t^{\odot s_{1}}(z) \odot \big[
   T(z) \big(a_{1} z \partial_{z} \mu^{[s_{2}]}(z) + a_{0}
   \mu^{[s_{2}]}(z)\big)\big] = \\
   \begin{cases}
      c
          n^{s \alpha - \frac{1}{2}} + O\big(n^{s \alpha -
          1}\big), &  s_{2} = 0 \\
          \binom{s}{s_{1}} \frac{s_{2}! \sigma^{s_{2}-1}
          \tau}{2^{(s_{2}+1)/{2}} 
          \sqrt{\pi}} \prod_{j=1}^{s_{2}}\frac{\Gamma(j
          \alpha')}{\Gamma(j \alpha' + \frac{1}{2})} n^{s \alpha +
          \frac{s_{2}-1}{2}}
          + O\big(n^{s \alpha + \frac{s_{2}-1}{2}- q}\big), &
          1 \le s_{2} < s.
   \end{cases}
\end{multline*}

Thus under the  assumptions given above, the dominant contribution
to~$r^{[s]}(z)$ is obtained when $s_{2} = s-1$ and  $s_{1}=1$,
giving the expansion 
\begin{equation*}
   \rho^{n} [z^{n}] r^{[s]}(z) = 
   \frac{s! \sigma^{s-2} \tau}{2^{{s}/{2}} \sqrt{\pi}}
   \prod_{j=1}^{s-1} \frac{\Gamma(j \alpha')}{\Gamma(j \alpha' +
   \frac{1}{2})} n^{s\alpha'-1}
   + O\big(n^{s \alpha' - 1 - q}\big),
\end{equation*}
which in turn yields the following singular expansion for 
$r^{[s]}(z)$:
\begin{equation}
   r^{[s]}(z) = \frac{s! \sigma^{s-2} \tau}{2^{{s}/{2}}
   \sqrt{\pi}} 
   \gamma_s Z^{-s \alpha'} 
   + O\big(|Z|^{-s \alpha'+ q}\big).
\end{equation}
Immediately from~\eqref{eq:26} and~\eqref{eq:10} follow the expansions
\begin{equation*}
   g^{[s]}(z) = \frac{s! \sigma^{s-1} \tau}{2^{({s+1})/{2}} \sqrt{\pi} \rho}
   \gamma_s Z^{-s \alpha' - \frac{1}{2}} 
   + O\big(|Z|^{-s \alpha' - \frac12 + q}\big)
\end{equation*}
and
\begin{equation*}
  \frac{g^{[s]}(z)}{T(z)} = \frac{s! \sigma^{s-1}}{2^{({s+1})/{2}}
   \sqrt{\pi} \rho} 
   \gamma_s Z^{-s \alpha' - \frac{1}{2}} 
   + O\big(|Z|^{-s \alpha' - \frac12 + q}\big).
\end{equation*}
Integrating leads to
\begin{equation*}
   \int_{0}^{z} \frac{g^{[s]}(t)}{T(t)} dt = \frac{s!
   \sigma^{s-1}}{2^{({s+1})/{2}}  
   \sqrt{\pi} (s \alpha' - \frac{1}{2})} 
   \gamma_s Z^{-s \alpha' + \frac{1}{2}} 
   + O\big(|Z|^{-s \alpha'+ \frac12 + q}\big).
\end{equation*}
Using~\eqref{eqno4} and~\eqref{eq:7}, 
we obtain~\eqref{eqno6}
and Theorem~\ref{theo1} is proved for $\alpha > 0$. 

The case $\alpha=0$
has already been
proved 
in \cite{panholzer:2003}, where the distribution has been
characterized by its moments.  Therefore we describe only very 
briefly how to obtain this result with the present approach.

One need only show by induction the singular behavior
\begin{multline}
   \label{eqno7}
   \mu^{[s]}(z) \\=
   \begin{cases}
   \frac{\tau}{2} \ln Z^{-1} + \A + O\big(|Z^{{1}/{2}} \log Z^{-1}|
   \big),     & s=1, \\
   \sqrt{2} \sigma \tau Z^{-{1}/{2}} 
   + O\big(|\log Z^{-1}|^2\big), & s=2, \\
   \frac{\tau}{\sqrt{\pi}} \sigma^{s-1} 2^{({s-1})/{2}}
   \Gamma(\frac{s}{2}+1) 
   \Gamma(\frac{s-1}{2}) Z^{-({s-1})/{2}} +
   O\big(|Z^{-\frac{s}{2}+1} \log Z^{-1}|\big), & s \ge 3.
   \end{cases}
\end{multline}
The desired result then follows by applying singularity analysis and
the duplication formula for the $\Gamma$-function.

To begin the proof of~\eqref{eqno7}, first we remark that for $s=1$
one proceeds as in Section~\ref{seco1}
and gets the full expansion
 \begin{equation*}
    \mu^{[1]}(z) = \frac{\tau}2 \ln Z^{-1} +  \A + (Z^{{1}/{2}} \log
          Z^{-1}) \A    + Z^{{1}/{2}} \A +  (Z \log Z^{-1}) \A,
 \end{equation*}
which  of course 
gives~\eqref{eqno7} in that case.
Assuming that~\eqref{eqno7} holds for $1 \le s_{2} < s$ with a given
$s \ge 2$,  we have the singular expansion
\begin{equation*}
  \begin{split}
   \partial_{z} \mu^{[s_{2}]}(z)  = 
   \frac\tau{\sqrt\pi\rho}{\sigma^{s_{2}-1} 2^{({s_{2}-1})/{2}}
   \Gamma\left(\frac{s_{2}}{2}+1\right) \Gamma\left(\frac{s_{2} +
   1}{2}\right)} 
   Z^{-({s_{2}+1})/{2}} \\
   {}+ O\big(|Z^{-\frac{s_2}{2}} \log Z^{-1}|\big).
  \end{split}
\end{equation*}

Under the restrictions $s_{1}+s_{2}=s$ and $s_{2}<s$, we obtain via
singularity analysis the expansions
\begin{multline*}
   \rho^{n} [z^{n}] \left\{\binom{s}{s_{1}} t^{\odot s_{1}}(z) \odot \big[T(z) 
   \big(a_{1} z \partial_{z} \mu^{[s_{2}]}(z) + a_{0}
   \mu^{[s_{2}]}(z)\big) \big]\right\} \\ =
      \binom{s}{s_{1}} \frac{\tau}{\sqrt\pi} \sigma^{s_{2}-1}
      2^{\frac{s_{2}-1}{2}} \Gamma\left(\frac{s_{2}}{2}+1\right)
      n^{\frac{s_{2}-1}{2}}  + O\big(n^{\frac{s_{2}}{2}-1} \log
      n\big)
\end{multline*}
for $1 \le s_{2} < s$.  [For $s_2 = 0$, an expansion is already
available at~\eqref{eq:33}.]
Under the given restrictions, the dominant contribution to
$r^{[s]}(z)$  is obtained when  $s_{2}=s-1$ 
and we obtain the following singular behavior 
of $r^{[s]}(z)$:
\begin{equation*}
   r^{[s]}(z) = \frac\tau{\sqrt\pi} {s \sigma^{s-2}
   2^{\frac{s}{2}-1} 
   \Gamma\left(\frac{s+1}{2}\right) \Gamma \left( \frac{s}{2} \right)}
   Z^{-{s}/{2}} +  O\big(|Z^{-\frac{s-1}{2}} \log Z^{-1}|\big).
\end{equation*}

We proceed with
\begin{equation*}
   \frac{g^{[s]}(z)}{T(z)} = \frac{1}{\rho\sqrt\pi} {s \sigma^{s-1}
   2^{\frac{s-3}{2}} \Gamma\left(\frac{s+1}{2}\right) 
   \Gamma\left(\frac{s}{2}\right)} Z^{-\frac{s+1}{2}}
   + O\big(|Z^{-{s}/{2}} \log Z^{-1}|\big),
\end{equation*}
and, due to~\eqref{eqno4} and~\eqref{eq:7}, integrating
gives~\eqref{eqno7} for $s \geq 2$ and completes the proof.
\end{proof}
\bibliographystyle{habbrv}
\bibliography{master}

\end{document}